\documentclass{amsart}
\usepackage{amssymb,enumerate,verbatim,xspace}
\usepackage[latin1]{inputenc}
\newcommand{\sscl}[1]{\cl(#1)}
\def\C{\mathcal C}
\def\M{\mathfrak M}
\def\N{\mathfrak N}
\def\cl{\mathrm{cl}}
\def\LL{\mathcal L}
\def\FF{\mathbb F}
\def\QQ{\mathbb Q}
\def\acl{\mathrm{acl}}
\def\gcl{\mathrm{gcl}}
\def\p{\varphi}
\def\tr{\mathrm{tr.deg}}
\def\cb{\mathrm{Cb}}
\def\ldim{\mathrm{lin.dim}_{\mathbb F_p}}
\def\Th{\mathrm{Th}}
\def\tp{\mathrm{tp}}
\def\U{\,\ddot{\!U}}
\def\mdim{\mathrm{lin.dim}_{\mathbb Q}}
\def\lin{\mathrm{lin.dim}}

\def\qftp{\mathrm{qftp}}
\def\Ind#1#2{#1\setbox0=\hbox{$#1x$}\kern\wd0\hbox to 0pt{\hss$#1\mid$\hss}
\lower.9\ht0\hbox to 0pt{\hss$#1\smile$\hss}\kern\wd0}
\def\Notind#1#2{#1\setbox0=\hbox{$#1x$}\kern\wd0\hbox to 0pt{\mathchardef
\nn="3236\hss$#1\nn$\kern1.4\wd0\hss}\hbox to 0pt{\hss$#1\mid$\hss}\lower.9\ht0
\hbox to 0pt{\hss$#1\smile$\hss}\kern\wd0}
\def\ind{\mathop{\mathpalette\Ind{}}}
\def\nind{\mathop{\mathpalette\Notind{}}}

\begin{document}

\title[Fields and Fusions]{Fields and Fusions\\Hrushovski constructions and their definable groups}

\author{Frank O Wagner}
\thanks{Membre junior de l'Institut Universitaire de France.\\
Partially supported by ANR-09-BLAN-0047 MODIG}
\date{\today}
\begin{abstract}An overview is given of the various expansions of fields and fusions of strongly minimal sets obtained by means of Hrushovski's amalgamation method, as well as a characterization of the groups definable in these structures.
\end{abstract}
\maketitle

\section{Introduction}
In 1986 Ehud Hrushovski invented a new method to obtain new stable structures from the class of their finitely generated substructures, via an adaptation of Fra\"\i ss\'e's construction of a universial homogeneous countable relational structure from its finite substructures. In particular, he constructed an $\aleph_0$-categorical stable complete pseudoplane (refuting a conjecture of Lachlan), a strongly minimal set with an exotic geometry which is not disintegrated, but does not interpret any group (refuting a conjecture of Zilber), and the fusion of two strongly minimal sets in disjoint languages in a third one (proving the non-existence of a maximal strongly minimal set).

His method was taken up by a number of people who adapted the technique to construct various exotic objects. Most recently, after preliminary work by Poizat \cite{Po99,Po01}, Baudisch, Mart\'\i n Pizarro and Ziegler achieved the fusion of two strongly minimal sets over a common $\FF_p$-vector space \cite{BMPZredfus} and a Morley rank $2$ expansion of an algebraically closed field of positive characteristic by a predicate for an additive non-algebraic subgroup \cite{BMPZred}, and Baudisch, Mart\'\i n Pizarro, Hils and the author constructed a Morley rank $2$ expansion of an algebraically closed field of characteristic zero by a predicate for a multiplicative non-algebraic subgroup \cite{BHMW}, a so-called {\em bad field}. I shall describe the basic construction and give some details on how the coloured fields are obtained.

A recurring question in the model-theoretic analysis of a structure is that of characterizing its definable groups. For the initial constructions, Hrushovski has answered it almost completely: The $\aleph_0$-categorical pseudoplane and the new strongly minimal set do not interpret any group, and a group definable in the fusion of two strongly minimal sets in disjoint languages is isogenous to a product of two groups definable in either set. In \cite{BMPW09} Blossier, Mart\'\i n Pizarro and the author introduced a geometric property for a theory relative to a reduct, relative CM-triviality, which holds in particular for many structures obtained by Hrushovski amalgamation and which yields a description of the definable groups. It follows in particular that a simple group definable in a fusion of two strongly minimal sets over an $\aleph_0$-categorical reduct is definable in one of the sets, and a simple group definable in a coloured field is algebraic. More generally, modulo a central subgroup,\begin{enumerate}
\item a group definable in a coloured field is an extension of a subgroup of some Cartesian power of the colour subgroup by an algebraic group; 
\item a group definable in the fusion of strongly minimal $T_1$ and $T_2$ over a common $\FF_p$-vector space is an extension of a definable $\FF_p$-vector space by a product of a $T_1$-definable and a $T_2$-definable group.\end{enumerate}

Most of this survey was presented at the Kirishima Model Theory Meeting 1-5 March 2010, one of four Camp-Style Seminars funded by RIMS of Kyoto University. I should like to thank the organizers Koichiro Ikeda (Hosei University), Masanori Itai (Tokai University, chair), Hirotaka Kikyo (Kobe University) and Akito Tsuboi (University of Tsukuba) for a very enjoyable meeting.

\section{The basic construction}
\subsection{Fra\"\i ss\'e's original construction}
Let $\C$ be a class of finite structures in a finite relational language, closed under substructures, and with the {\em amalgamation property} AP (where we allow $A=\emptyset$):
\begin{quote} For all injective $\sigma_i:A\rightarrow B_i$ in $\C$ for $i=1,2$ there are injective $\rho_i:B_i\rightarrow D\in\C$ with $\rho_1\sigma_1=\rho_2\sigma_2$.\end{quote}
Then there is a unique countable structure $\M$ such that:
\begin{quote} For all finite $A\subset\M$ and $A\subset B\in\C$ there is an embedding $B\rightarrow\M$ which is the identity on $A$.\end{quote}
The proof is by successive amalgamation over all possible situations, using AP. Note that $\C$ is countable.

We call $\M$ the {\em generic} model; it is ultrahomogeneous, and hence $\aleph_0$-categorical (since the language is finite).

Fra\"\i ss\'e's aim was to construct $\aleph_0$-categorical structures. If we do not mind losing $\aleph_0$-categoricity, we can drop various conditions:\begin{itemize}
\item We can work in an arbitrary language (but we shall stick to countable ones for simplicity).
\item We can either work with a class of finitely generated structures (in which case we should verify that {\em all} substructures of a finitely generated structures are finitely generated),
\item or with a class of algebraic closures of finitely generated structures (typically, algebraically closed fields of finite transcendence degree),
\item or with a class of countable structures, of size at most $\aleph_1$ (this will usually happen in the strictly stable case).\end{itemize}
In any case our class should be axiomatizable (apart from the cardinality restriction); since it is closed under substructure, the axiomatization will be universal.

Assuming $\aleph_1=2^{\aleph_0}$, let $\C^+$ be the class of structures of size at most $\aleph_1$ whose countable substructures are in $\C$. A structure in $\C^+$ can then be written as an increasing union of structures in $\C$. Using AP for $\C$ iteratively, we have:
\begin{quote} For all injective $\sigma_i:A\rightarrow B_i$ with $A,B_2\in\C$ and $B_1\in\C^+$ there are injective $\rho_i:B_i\rightarrow D\in\C^+$ with $\rho_1\sigma_1=\rho_2\sigma_2$.\end{quote}
Then there is a unique structure $\M$ of size $\aleph_1$ such that:
\begin{quote} For all $A\subset\M$ and $A\subset B\in\C$ there is an embedding $B\rightarrow\M$ which is the identity on $A$.\end{quote}
We still call $\M$ the generic model; it is ultrahomogeneous for countable substructures.
(If $\aleph_1\not= 2^{\aleph_0}$, work with structures of size at most $\aleph_\alpha$, where $\alpha$ is minimal with $\aleph_\alpha^+\ge|\C|$.)

\subsection{Strong embeddings} Rather than considering all inclusions $A\subset B\in\C$, we only consider {\em certain} inclusions $A\le B$, which we call {\em strong}. We require $\le$ to be transitive and preserved under intersections. We only demand AP for strong embeddings, and obtain a generic structure $\M$ such that the {\em richness} condition holds:
\begin{quote} For any $A\le\M$ and $A\le B\in\C$ there is a {\em strong} embedding $B\rightarrow\M$ which is the identity on $A$.\end{quote}
For $A\subset\M$ with $A\in\C$ we define the {\em closure} $\cl_\M(A)$ to be the smallest $B\le\M$ containing $A$ (note that the closure exists since strongness is preserved under intersections, but it might be infinite). We should choose $\C$ so that the closure of a set in $\C$ is again in $\C$.
Then $\M$ is ultrahomogeneous for {\em strong} subsets in $\C$. Note that if $\M$ is sufficiently saturated, then by uniqueness $\cl_\M(A)$ is contained in the model-theoretic algebraic closure of $A$ in the sense of $\Th(\M)$.

In order to axiomatize $\Th(\M)$, we need to express $A\le\M$. If this can be done by a first-order formula, and if there is a bound on the number of possible closures (and hence types) of a finite subset of $\M$, the generic model is $\aleph_0$-categorical.

However, even in a finite relational language for finite substructures, closedness need not be a definable property. We shall require it to be {\em type-definable}, and we need {\em approximate definability} of richness:
\begin{quote} If $A\le B\in\C$ and $A$ is sufficiently strong in $\M$, then there is an embedding of $B$ into $\M$ over $A$ whose image has a pre-described leved of strength.\end{quote}
This yields richness, and thus homogeneity for countable strong subsets of a $\aleph_1$-saturated model.

The axiomatization usually requires the generic model to be sufficiently saturated, since then a structure will be generic for $\C$ if and only if it is a $\aleph_1$-saturated model of $\Th(\M)$.
In order to simplify the exposition, we shall from now on assume the generic model to be $\aleph_1$-saturated. However, there are interesting examples of non-saturated generic structures \cite{ShSp,BS,L,Ik,I}.

\subsection{Predimension} In order to define a strong embedding relation related to a rudimentary notion of dimension, we consider a function $\delta$ from the set of all finitely generated structures in $\C$ to the non-negative reals, satisfying $\delta(\emptyset)=0$ and
$$\mbox{(Submodularity)}\qquad
\delta(A)+\delta(B)\ge\delta(AB)+\delta(A\cap B).$$
More generally, we define a relative predimension $\delta(A/C)$ for a set $A$ finitely generated over some parameters $C$ such that $\delta(A/C)=\delta(AC)-\delta(C)$ if $C$ is finite, and
$$\delta(A/C)+\delta(B/C)\ge\delta(AB/C)+\delta((A\cap B)/C).$$
Define
$$A\le B\quad\Leftrightarrow\quad\delta(B'/A)\ge0\mbox{ for all finitely generated $B'\subseteq B$}.$$
Submodularity the ensures that $\le$ is transitive and closed under intersections. By considering a suitable limit, we can define a {\em dimension} $d_\M(A/C)=\delta(\cl_\M(AC)/C)$, which will be submodular for strong subsets.

Note that if $B\le\M$ and $\delta(\bar a/B)=0$, then $\bar aB\le\M$.

Usually, $\delta(\bar a/A)\le r$ is a closed condition in the following strong sense: Given $A\le\M$ and  $r'\ge 0$ there is a collection $\Phi_{r'}(A)$ of (quantifier-free) $\LL(A)$-formulas with $\delta(\bar a'/A)\le r'$ for any $\p\in\Phi_{r'}(A)$ and $\bar a'\models\p$, and such that
$$\delta(\bar a/A)=\inf\{r':\tp(\bar a/A)\cap\Phi_{r'}(A)\not=\emptyset\}.$$
This is very useful to axiomatize approximate richness, but I do not know whether it is necessary. Note that this implies the existence of a collection $\Psi_{r'}(A)$ of existential formulas with $d_\M(\bar a'/A)\le r'$ for any $\p\in\Psi_{r'}(A)$ and $\bar a'\models\p$, such that
$$d_\M(\bar a/A)=\inf\{r':\tp(\bar a/A)\cap\Psi_{r'}(A)\not=\emptyset\}$$
(quantify existentially over the elements in the closure).

Given $A\le\M\prec\N$ (with $\N$ saturated) and a type $p(\bar x)\in S(A)$, put
$$d(p)=\sup\{d_\N(\bar a/M):\bar a\models p\}.$$
Then for all $r<d(p)$ there is $\bar a\models p$ in some elementary extension of $\M$ with $\models\neg\p(\bar a)$ for all possible formulas $\p\in\Psi_r(M)$. By compactness there is a realization $\bar a\models p$ not realizing any formula in $\Psi_r(M)$ for $r<d(p)$, so $d_\N(\bar a/M)=r$ and the supremum is attained. Clearly $d(\tp(\bar a/A))\le d_\M(\bar a/A)$, but I do not see a reason why it might not sometimes be strictly smaller (unless $A$ is itself a generic model).

Now suppose $B\le C\le\M$ and $\delta(\bar a/B)=\delta(\bar a/C)$, with $\langle\bar a B\rangle\le\M$. Then $\langle\bar aC\rangle$ is also strong, since for any $C'\supseteq C$
$$\delta(\bar a/C')\ge\delta(\bar a/B)=\delta(\bar a/C).$$
So $d(\bar a/B)=\delta(\bar a/B)=\delta(\bar a/C)=d(\bar a/C)$.
Put $B'=\langle\bar a B\rangle\cap C$, a strong subset. Then
$$\delta(\bar a/B)+\delta(C/B)\ge\delta(\bar aC/B)+\delta(B'/B)\\
=\delta(\bar a/C)+\delta(C/B)+\delta(B'/B)$$
by submodularity, whence
$$\delta(\bar a/B)\ge\delta(\bar a/C)+\delta(B'/B)=\delta(\bar a/B)+\delta(B'/B).$$
(If $\delta(C/B)$ is infinite, approximate by closures of finitely generated subsets over $B$.)
Since $B\le\M$ implies $\delta(B'/B)\ge0$ we get $\delta(B'/B)=0$.
If $B'=B$ we call $\tp(\bar a/C)$ a {\em free} extension of $p$ to $\M$, and we say that $\langle\bar a B\rangle$ and $C$ are {\em freely amalgamated} over $B$.

If the class $\C$ is closed under free amalgamation, we call its generic structure {\em free}.

So far the trivial predimension $\delta\equiv 0$ has not ben excluded. (For the trivial predimension any containment of sets is strong, so we are back to the Fra\"\i ss\'e method.) In order to obtain {\em stable} structures, we require the free extension to $\M$ of a type over an algebraically closed set to be unique, in which case it will be the non-forking extension. We can now count types: For any $\aleph_1$-saturated model $\M$ and any $p\in S(M)$ we chose a formula $\p_r\in p\cap\Psi_r(M)$ for all rational $r\ge d(p)$. Let $A_0\le M$ be the algebraic closure of the parameters used in these formulas, $\bar a\models p$ and $A=\langle\bar a A_0\rangle\cap M$. Then $d(p{\restriction_A})=d(p)$ and $\langle\bar a A\rangle\cap M=A$, so $p$ is the unique free extension of $p{\restriction_A}$. But this means that there are at most $|\M|^{\aleph_0}$ types over $\M$, and $\M$ is stable. 

For generic models $\M\prec\N$ we have $\bar a\ind_\M\N$ if and only if $d(\bar a/\M)=d(\bar a/\N)$ and $\langle\bar a\N\rangle$ is strong. In other words, an extension of the same dimension can only fork if $\langle\M\bar a\rangle\cap\N\supset\M$, i.e.\ some non-algebraic elements in its closure become algebraic. In particular, a type of dimension $0$ can only fork by algebraicizing elements in its closure.

\subsection{Example (Ab initio)}
In a relational language $\LL$, choose a {\em weight} $\alpha_R>0$ for every relation $R$, and define a {\em predimension} on finite structures:
$$\delta(A)=|A|-\sum_{R\in\LL}\alpha_R|R(A)|,\qquad\mbox{as well as}$$
$$\delta(A/B)=\delta(AB)-\delta(B)=|A\setminus B|-\mbox{weights of the new relations};$$
the latter makes sense even if $B$ is infinite.

Let $\C$ be the universal class of all finite $\LL$-structures whose substructures have non-negative predimension. It is closed under free amalgamation (an amalgam of $A$ and $B$ over their intersection $A\cap B$ is {\em free} if any relation on $A\cup B$ either lives on $A$ or on $B$), and thus has AP. Closedness is type-definable (this uses $\delta\ge0$), richness is approximately definable, and $\aleph_1$-saturated models are rich. The generic model $\M$ is $\omega$-stable if all the $\alpha_R$ are multiples of a common rational fraction, and stable otherwise; two strong subsets are independent in the forking sense if and only if they are freely amalgamated over their intersection and the amalgam is strong in $\M$. It is now not to difficult to see that $\M$ weak elimination of imaginaries.

Usually, one wants to consider subclasses of $\C$ with AP, in order to obtain $\aleph_0$-categorical or strongly minimal strutures. 

\subsection{Morley rank}
Suppose the range of $\delta$ is closed and discrete. Then every type over $\M$ has a restriction to a finite subset of which it is the free extension. So $\Th(\M)$ is superstable; it will be $\omega$-stable if multiplicities are finite.

For any strong $A\le\M$ and tuple $\bar a\in\M$ there is a (finite) tuple $\bar b$ such that $$\cl_\M(\bar aA)=\langle A\bar a\bar b\rangle.$$
As $\cl_M(\bar aA)\subseteq\acl(\bar aA)$, the type $\tp(\bar a/A)$ determines $\tp(\bar a\bar b/A)$. Conversely, the quantifier-free type $\qftp(\bar a\bar b/A)$ determines $\tp(\bar a/A)$ by homogeneity.

Suppose $\delta(\bar a\bar b/A)=d_\M(\bar a/A)=r$ and $\p_r(x)\in\qftp(\bar a\bar b/A)$ is such that $\p_r$ isolates $\tp(\bar a\bar b/A)$ among the types of predimension at least $r$. Then $\exists\bar y\,\p_r(\bar x,\bar y)$ isolates $\tp(\bar a/A)$ among types of dimension at least $r$. If enough such formulas exist (for instance, in the {\em ab initio} case), $\Th(\M)$ is $\omega$-stable.

A proper strong extension $A<B\in\C$ is {\em minimal} if $A\le A'\le B$ implies $A'=A$ or $A'=B$. Equivalently, $\delta(B/A')<0$ for all $A\subset A'\subset B$. The extension is {\em pre-algebraic} if $\delta(B/A)=0$.

For a minimal pre-algebraic extension $A<B$ let $A_0\le A$ be the closure of the canonical base
$\cb(B\setminus A/A)$. This is the unique minimal strong subset of $A$ over which $B\setminus A$ is pre-algebraic (and in fact minimal). We call $A_0<(B\setminus A)\cup A_0$ {\em bi-minimal pre-algebraic}. Note that a minimal pre-algebraic extension can only fork by becoming algebraic; it must thus have Lascar rank $1$.

Let $A\le B\in\C$ be pre-algebraic. If $B_0\le B$ has predimension $0$ over $A$ and minimal Lascar rank possible, it must be minimal. Hence $U(B_0/A)=1$; since $U(B/AB_0)<U(B/A)$, we see inductively that $\tp(B/A)$ has finite Lascar rank.

Since a forking extension of a general type can only have the same dimension if some tuple in its closure of dimension $0$ becomes algebraic, well-foundedness of dimension implies that every type has Lascar rank $<\omega^2$.

\subsection{Geometry and the collapse}
Suppose that $\delta$ takes integer values, and the maximal dimension of a point is $1$. For a single point $a$ and a set $B\le\M$ there are two possibilities:\begin{itemize}
\item $d_\M(a/B)=\delta(a/B)=1$. Then $aB\le\M$, so this determines a unique type, the {\em generic} type.
\item $d_\M(a/B)=0$. So $a$ is in the {\em geometric} closure $\gcl(B)$.\end{itemize}

Clearly $\gcl$ is increasing and idempotent, hence a closure operator, which in addition satisfies the exchange rule:
\begin{quote} If $a\in\gcl(Bc)\setminus\gcl(B)$, then $c\in\gcl(Ba)$.\end{quote}
Note that the generic structure has rank (Lascar rank, and Morley rank if it is $\omega$-stable) at most $\omega$, since a forking extension of the generic type has dimension $0$, and thus finite rank. More generally, if $\delta$ takes integer values and a single point has maximal dimension $d$, then the rank is bounded by $\omega\cdot d$.

We should like to restrict the class $\C$ so that $\gcl$ becomes algebraic closure, thus yielding a {\em strongly minimal} set (every definable subset is uniformly finite or co-finite). For that, we have to bound the number of possible realisations of any $a\in\gcl(B)$.
Clearly it is sufficient to bound the number of realisations for each bi-minimal pre-algebraic extension.
So let $\mu$ be a function from the set of isomorphism types of bi-minimal pre-algebraic extensions to the integers, and let $\C^\mu$ be the class of $A\in\C$ such that for any such extension $A_0<B$ with $A_0\le A$ there are at most $\mu(A_0,B)$ independent strong copies of $B$ in $A$ over $A_0$. In order to show AP we usually verify that $\C^\mu$ has {\em thrifty amalgamation}:
\begin{quote}
If $A\le B\in\C^\mu$ is minimal and $A\le M\in\C^\mu$, then either the free amalgam of $A$ and $\M$ over $B$ is still in $\C^\mu$, or $B$ embeds closedly into $M$ over $A$.\end{quote}
Then a generic model for $\C^\mu$ exists (even as a strong substructure of a generic model for $\C$).

However, in order to axiomatize the class $\C^\mu$, we have to restrict the number of independent realizations of bi-minimal pre-algebraic extensions {\em uniformly} and {\em definably}. In other words, for every such extension $A_0<B$ and for $n=\mu(A_0,B)$ we choose a (quantifier-free) formula $\psi(X,Y_0,\ldots,Y_n)\in\tp(A_0,B_0,\ldots,B_n)$, where $B_0,\ldots,B_n$ are $n+1$ independent realizations of $\tp(B/A_0)$, and we consider the class $\C_\mu$ given by the (universal) axioms
$$\forall X\neg\exists Y_0\ldots Y_n\,\psi(X,Y_0,\ldots,Y_n).$$
Suppose $\C_\mu$ is non-empty and has thrifty amalgamation. It will then inherit approximate definability of richness from $\C$ provided we can definably check that minimal extensions are in $\C_\mu$ and not just in $\C$. That is, given an extension $A<B$ with $A\in\C_\mu$ there should only be finitely many axioms $\forall X\neg\exists Y_0\ldots Y_n\,\psi(X,Y_0,\ldots,Y_n)$ which $B$ could possibly violate. If $\mu$ is finite-to-one, this is ensured by dimension considerations. So the generic structure for $\C_\mu$ is strongly minimal. It is called the {\em collapse} of the generic structure for $\C$.

\section{Fields and Fusion}
\subsection{Fusion}
\begin{itemize}
\item Two strongly minimal sets with the {\em definable multiplicity property} DMP in disjoint languages can be amalgamated freely with predimension
$$\delta(A/B)=RM_1(A/B)+RM_2(A/B)-|A\setminus B|$$
and collapsed to a strongly minimal set \cite{Hr93} (see also \cite{BMPZfus}).
\item Two theories of finite and definable Morley rank with DMP in disjoint languages can be amalgamated freely with predimension
$$\delta(A/B)=n_1\cdot RM_1(A/B)+n_2\cdot RM_2(A/B)-n\cdot|A\setminus B|$$
where $n_1\cdot RM(T_1)=n_2\cdot RM(T_2)=n$, and collapsed to a structure of Morley rank $n$ \cite{Zi}.
\item Two strongly minimal sets with DMP with a common $\aleph_0$-categorical reduct, one preserving multiplicities, can be amalgamated freely with predimension
$$\delta(A/B)=RM_1(A/B)+RM_2(A/B)-RM_0(A/B)$$
and collapsed to a strongly minimal set \cite{BMPZredfus} (see also \cite{HH} for partial results, and \cite{Hi} for the extension from a common vector space reduct to a common $\aleph_0$-categorical reduct).
\item {\em Conjecture.} Two theories of finite and definable Morley rank with DMP with a common $\aleph_0$-categorical reduct, one preserving multiplicities, can be amalgamated freely with predimension
$$\delta(A/B)=n_1\cdot RM_1(A/B)+n_2\cdot RM_2(A/B)-n_0\cdot RM_0(A\setminus B)$$
where $n_0\cdot RM(T_0)=n_1\cdot RM(T_1)=n_2\cdot RM(T_2)=n$, and collapsed to a structure of Morley rank $n$.
\end{itemize}
Note that in a strongly minimal set, Morley rank is always definable. More generally, definability of Morley rank yields approximate definability of closedness. We also need definability of rank, together with the DMP, for the collapse to be definable.

In order to ensure submodularity of the predimension, the negative part of the predimension should be {\em modular}, i.e.\ equality should hold. This is trivial for cardinality, and implied by $\aleph_0$-categoricity. However, $\aleph_0$-categoricity of the common reduct is not {\em necessary} for the free or collapsed fusion to exist (see the remarks on the green field below).

For the collapse over equality, or over a disintegrated reduct, it is relatively easy to find a suitable $\mu$-function. The essential ingredient is that any type of finite rank over a free amalgam of $B$ and $C$ over $A$ must be based on either $B$ or $C$. However, for the collapse over a common $\aleph_0$-categorical reduct, this is no longer true. One uses the fact that the reduct is essentially a finite cover of a vector space over a finite field, and then translates bi-minimal pre-algebraic types over a free amalgam so that they become based on the left or on the right. This is similar to the red field studied in \ref{redfield}.

A {\em coloured field} is a field expanded by a unary predicate, called its colour.\begin{itemize}
\item The colour {\em black} distinuishes an algebraically independent subset.
\item The colour {\em red} distinuishes a proper non-trivial connected additive subgroup.
\item The colour {\em green} distinuishes a proper non-trivial connected multiplicative subgroup.\end{itemize}
Green fields of finite Morley rank (so-called {\em bad fields}) first came up in early work on the Cherlin-Zilber conjecture that groups of finite Morley rank are algebraic groups over an algebraically closed field.

I shall describe the construction of the red and of the green field as separate constructions.
However, we can also comprehend:\begin{itemize}
\item The black field as the fusion over equality of an algebraically closed field and the structure consisting of a pure set with an infinite co-infinite predicate. (Actually, rank considerations might be easier if we have infinitely many infinite disjoint predicates, and in the end throw all but one of them away.)
\item The red field as the fusion of an algebraically closed field of characteristic $p>0$ and an elementary abelian $p$-group with a predicate for an infinite subgroup of inifinite index, over the common reduct to the abelian (additive) $p$-group.
\item The green field as the fusion of an algebraically closed field with the theory of a divisible abelian group (of the right torsion) with a distinguished torsion-free subgroup, over the pure (multiplicative) group structure as common reduct. Note that this kind of fusion has not been done in general, as the common reduct is not $\aleph_0$-categorical. Its existence, and certainly its axiomatization, depends on particular algebraic properties; one would have to distill the precise conditions (normally implied by $\aleph_0$-categoricity of the reduct) needed to make the fusion work.\end{itemize}

\subsection{Red fields}\label{redfield}
A {\em red field} is an $\omega$-stable algebraically closed field $K$ with a predicate $R$ for a connected additive subgroup of comparable rank.
Note that in characteristic $0$ this gives rise to an infinite definable subfield $\{a\in K:aR\le R\}$, so the structure has rank at least $\omega$. Hence we assume that the characteristic is positive.

Let $\C$ be the class of finitely generated fields $k$ of characteristic $p>0$ with a predicate $R$ for an additive subgroup, the {\em red} points, such that for all finitely generated subfields $k'$
$$\delta(k')=2\,\tr(k')-\ldim(R(k'))\ge 0.$$
This condition is universal, since we have to say that $2n$ linearly independent red points do not lie in any variety of dimension $<n$.

For $k\le k'\in\C$ we define a predimension $\delta$ by
$$\delta(k'/k)=2\,\tr(k'/k)-\ldim(R(k')/R(k)).$$
Since $\C$ has free amalgamation, a generic model $\M$ exists; as richness is approximately definable, $\aleph_0$-saturated models of $\Th(\M)$ are rich.

For a point $a$ and a set $B$ in $\M$ there are three possibilities:\begin{enumerate}
\item $d_\M(a/B)=2$. Then $a$ is not red, and $aB$ is strong. $RM(a/B)=\omega\cdot2$.
\item $d_\M(a/B)=1$. There is a red point $a'$ interalgebraic with $a$ over $B$, and $a'B$ is strong. $\omega\cdot2>RM(a/B)\ge RM(a'/B)=\omega$.
\item $d_\M(a/B)=0$. Then either $a$ is algebraic over $B$, or pre-algebraic.\end{enumerate}

In order to collapse, we want to restrict the number of bi-minimal pre-algebraic extensions. A {\em code} is a formula $\p(\bar x,\bar y)$ with $n=|\bar x|$ such that\begin{enumerate}
\item For all $\bar b$ either $\p(\bar x,\bar b)$ is empty, or has Morley degree $1$. So $\p(\bar x,\bar b)$ determines a unique generic type $p_{\p(\bar x,\bar b)}$ (or is empty).
\item If $RM(\p(\bar x,\bar b)\cap\p(\bar x,\bar b'))=n/2$, then $b=b'$. In other words, $\bar b$ is the canonical base for $p_{\p(\bar x,\bar b)}$ and the extension $\bar b\le\bar a\bar b$ is bi-minimal.
\item $RM(\bar a/\bar b)=n/2$ and $\ldim(\bar a/\bar b)=n$ for generic $\bar a\models\p(\bar x,\bar b)$, and $2\,\tr(\bar a/U\bar b)<n-\ldim(U)$ for all non-trivial subspaces $U$ of $\langle\bar a\rangle$. Thus $\bar b\le\bar a\bar b$ is minimally pre-algebraic. Moreover, $\delta(\bar a'/B)<0$ for any $B\ni\bar b$ and non-generic $\bar a'\notin\acl(B)$ realizing $\p(\bar x,\bar b)$.
\item For any $H\in\mathrm{GL}_n(\FF_p)$, $\bar m$ and $\bar b$ there is $\bar b'$ with $\p(H\bar x+\bar m,\bar b)\equiv\p(\bar x,\bar b')$. Hence affine transformations preserve the code.
\item If $\p(\bar x,\bar b)$ is disintegrated for some $\bar b$, it is disintegrated (or empty) for all $\bar b$. So $\p$ fixes the type of the extension: disintegrated, or generic in a group coset (minimal pre-algebraic types are locally modular).
\end{enumerate}
By definability of Morley rank and multiplicity in algebraically closed fields these are definable properties (defining a coset of a group is definable by a Lemma of Ziegler \cite{ZiNote} which says that $p$ is a generic type of a coset of a subgroup of $(K^+)^n$ if and only if for independent realizations $a,b\models p$ the sum $a+b$ is independent of $a$). Enumerating all isomorphism types of bi-minimal pe-algebraic extensions, it is easy to find a set $\mathcal S$ of codes such that every minimal pre-algebraic extension is coded by a unique $\p\in\mathcal S$.

For a code $\p$ and some $\bar b$ consider a Morley sequence $(\bar a_0,\bar a_1,\ldots,\bar a_k,f)$ for $p_{\p(\bar x,\bar b)}$, and put $\bar e_i=\bar a_i-\bar f$.
We can then find a formula $\psi_\p^k\in\tp(\bar e_0,\ldots,\bar e_k)$ such that\begin{enumerate}
\item Any realization $(\bar e'_0,\ldots,\bar e_k')$ of $\psi_\p^k$ is $\mathbb F_p$-linearly independent, and $\models\p(\bar e'_i,\bar b')$ for some unique $\bar b'$ definable over sufficiently large finite subsets of the $\bar e'_i$, the {\em canonical parameter} of the sequence $\bar e'_0,\ldots,\bar e_k'$.
\item $\psi_\p^k$ is invariant under the finite group of {\em derivations} generated by
$$\partial_i:\bar x_j\mapsto\left\{\begin{array}{ll}\bar x_j-\bar x_i&\mbox{ if }j\not=i\\
-\bar x_i&\mbox{ if }j=i\end{array}\right.\mbox{ for }0\le i\le k.$$
\item Some condition ensuring dependence of affine combinations, and invariance under the stabiliser of the group for coset codes.\end{enumerate}

Given a code $\p$ and natural numbers $m,n$, there is some $\lambda$ such that for every $M\le N\in\C$ and realization $\bar e_0,\ldots,\bar e_\lambda\models\psi_\p^\lambda$ in $N$ with canonical parameter $\bar b$, either\begin{itemize}
\item the canonical parameter for some derived sequence lies in $M$, or
\item for every $A\subset N$ of size $m$ the sequence $(\bar e_0,\ldots,\bar e_\lambda)$ contains a Morley subsequence in $p_{\p(\bar x,\bar b)}$ over $MA$ of length $n$.\end{itemize}

Let $\mu$ be a sufficiently fast-growing finite-to-one function from $\mathcal S$ to $\omega$, and $\C_\mu$ the class of $A\in\C$ satisfying $\neg\exists\bar{\bar y}\,\psi_\p^{\mu(\p)}(\bar{\bar y})$ for all $\p\in\mathcal S$. The above lemma allows us to characterize when a minimal pre-algebraic extension of some $M\in\C_\mu$ is no longer in $\C_\mu$, and to prove thrifty amalgamation for $\C_\mu$. Hence there is a generic model $\M$, with $RM(\M)=2$ and $RM(R(\M))=1$. Alternatively to the standard axiomatization by richness, $\Th(\M)$ can be axiomatized inductively by\begin{itemize}
\item Finitely generated subfields are in $\C_\mu$.
\item ACF$_p$.
\item The extension of the model generated by a red generic realization of some code instance $\p(\bar x,\bar b)$ is not in $\C_\mu$.\end{itemize}
Since any complete theory of fields of finite Morley rank is $\aleph_1$-categorical, Lindstr\"om's theorem implies that $\Th(\M)$ is model-complete.

\subsection{Green fields}
A {\em green field} is an $\omega$-stable algebraically closed field $K$ with a predicate $\U$ for a connected multiplicative subgroup of comparable rank.
Note that in characteristic $p>0$ the existence of a green field of finite rank (which has $\tilde\FF_p$ as prime model) implies that there are only finitely many {\em $p$-Mersenne primes} $\frac{p^n-1}{p-1}$. Its existence is thus improbable; in any case it cannot be constructed as generic model by amalgamation methods \cite{Wa03}.

Let $\C$ be the class of finitely generated fields $k$ of characteristic $0$ with a predicate $\U$ for a torsion-free multiplicative subgroup, such that for all finitely generated subfields $k'$
$$\delta(k')=2\,\tr(k')-\mdim(\U(k'))\ge 0,$$
where the linear dimension is taken multiplicatively. Put
$$\delta(k'/k)=2\,\tr(k'/k)-\mdim(\U(k')/\U(k)).$$
While linear dimension over a finite field is definable, this is no longer true for dimension over $\QQ$, as there are infinitely many scalars (exponents).
Given a variety $V$ its {\em minimal torus} is the smallest torus containing $V$ in a single coset. We call a subvariety $W\subseteq V$
{\em cd-maximal} if its {\em codimension} $\mdim(W)-\tr(W)$ is strictly minimal among irreducible components of any $W'$ with $W\subset W'\subseteq V$.
Poizat used Zilber's Weak Intersections with Tori Theorem (weak CIT), a consequence of Ax' differential Schanuel conjecture, to show:\begin{quote}
For any uniform family $V_{\bar z}$ of varieties the set of minimal tori for its cd-maximal subvarieties is finite.\end{quote}
This specifies finitely many possibilities for $\QQ$-linear relations on a family of varieties which could render $\delta$ negative. Hence $\C$ is again universal and richness approximately axiomatizable.

However, we also have to worry about the size of $\C\,$! Suppose $\bar a$ is a generic point of some variety $V$ whose coordinates are green. Then for every $n<\omega$ there is a unique green $n$-th root of $\bar a$. Now if $\sqrt[n]{V}$ is not irreducible for infinitely many $n$, the type of $\bar a$ has to specify in which irreducible component of $\sqrt[n]{V}$ its green $n$-th root lies; this would yield $2^{\aleph_0}$ possible types. Fortunately, this does not happen \cite{Hi2}: For every $V$ there is $n$ (uniformly and definably in parameters) such that $\sqrt[k]{V}$ has at most $n$ irreducible components for all $k<\omega$.

Hence the generic model exists, and $\aleph_0$-saturated models of its theory are rich. It has Morley rank $\omega\cdot2$, and a generic green point has rank $\omega$. This also implies that the generic model has the DMP.

In order to collapse, we define codes similarly to the additive case. For Property (1), we have to ask for irreducibility of the set of $k$-th roots:\begin{itemize}
\item[$(1')$] For all $\bar b$ and all $k<\omega$ either $\p(\bar x^k,\bar b)$ is empty, or has Morley degree $1$.\end{itemize}
Property (2) remains unchanged. As we have to use the weak CIT in order to select a finite number of subspaces $U$ of $\langle\bar a\rangle$ we mention, property (3) becomes\begin{itemize}
\item[$(3')$] $RM(\bar a/\bar b)=n/2$ and $\mdim(\bar a/\bar b)=n$ for generic $\bar a\models\p(\bar x,\bar b)$, and for $i=2,\ldots,r$ and any $W$ irreducible component of $V\cap \bar a
T_i$ of maximal dimension, $\dim(T_i)>2\cdot\dim(W)$ if $V\cap \bar a
T_i$ is infinite.\end{itemize}
Since $\mathrm{GL}_n(\QQ)$ is infinite, we cannot encode invariance under the group of affine transformations in property (4) but have to treat it externally; we just demand invariance under multiplicative translation:
\begin{itemize}
\item[$(4')$] For any invertible $\bar m$ and $\bar b$ there is $\bar b'$ with $\p(\bar x\cdot\bar m,\bar b)\equiv\p(\bar x,\bar b')$.
\end{itemize}
Finally, since any algebraic subgroup of $(K^\times)^n$ is a torus and thus $\emptyset$-definable, all bi-minimal pre-algebraic extensions are disintegrated, eliminating the need for property (5).

Using the weak CIT we obtain:\begin{quote}
There exists a collection $\mathcal S$ of codes such that for every minimal pre-algebraic definable set $X$ there is a unique code $\p\in\mathcal S$ and finitely many tori $T$ such that $T\cap(X\times\p(\bar x,\bar b))$ projects generically onto $X$ and $\p(\bar x,\bar b)$ for some $\bar b$.\end{quote}
We call such a $T$ a {\em toric correspondence}. In particular, for any code $\p$ only finitely many tori can induce a toric correspondence between instances of $\p$.

For every code $\p$ and integer $k$ there is some formula $\psi_\p^k(\bar x_0,\ldots,\bar x_k)\in\tp(\bar e_0\cdot\bar f^{-1},\ldots,\bar e_k\cdot\bar f^{-1})$ for some Morley sequence $(\bar e_0,\ldots,\bar e_k,\bar f)$ in $p_{\p(\bar x,\bar b)}$ such that:\begin{itemize}
\item[$(1')$] Any realization $(\bar e'_0,\ldots,\bar e_k')$ of $\psi_\p^k$ is disjoint, and $\models\p(\bar e'_i,\bar b')$ for some unique $\bar b'$ definable over sufficiently large finite subsets of the $\bar e'_i$.
\item[$(2')$] If $\models\psi_\p^k(\bar e_0,\ldots,\bar e_k)$, then $\models\psi_\p^{k'}(\bar e_0,\ldots,\bar e_{k'})$ for each  $k'\leq k$, and $\psi_\p^k$ is invariant under derivations.
\item[$(3')$] Let $i\neq j$ and $(\bar e_0,\ldots,\bar e_k)$ realize $\psi$ with canonical parameter $\bar b$. If there is some toric correspondence $T$ on $\p$
and $\bar e_j'$ with $(\bar e_j,\bar e_j')\in T$, then $\bar e_i\nind_{\bar b}\bar e_j'\cdot\bar e_i^{-1}$ in case $\bar e_i$ is a generic realization of
$\p(\bar x,\bar b)$.\end{itemize}
Since $\QQ$ is infinite, we cannot demand $\QQ$-linear independence in property $(1')$ but merely disjointness. Similarly, in property $(3')$ we cannot check all linear combinations, but just the finitely many given by the toric correspondences.

This is sufficient to obtain the same counting lemma as before, characterize when minimal extension take us out of $\C_\mu$ (for some rapidly growing finite-to-one choice of $\mu$), prove thrifty amalgamation and axiomatize the theory of the generic model.
It has Morley rank $2$, and $\U^\M$ has Morley rank $1$. Moreover, there also is an alternative axiomatization analogous to the red case, which yields model-completeness.

\section{Definable groups}
\subsection{Relative CM-triviality} Let $T$ be a stable theory $T$ in a language $\LL$, and $T_0$ its reduct to a sublanguage $\LL_0$. (In the case of the fusion, we have two reducts and the definitions and results generalize to that case as well.) We assume that $T$ comes equipped with a finitary closure operator $\sscl.$ contained in the algebraic closure and satisfying\begin{enumerate} 
\item If  $A$ is algebraically closed and  $b\ind_Ac$, then $\sscl{Abc}\subseteq\acl_0(\sscl{Ab},\sscl{Ac})$.
\item If $\bar a\in\acl_0(A)$, then $\sscl{\acl(\bar a),A}\subseteq\acl_0(\acl(\bar a),\sscl{A})$.
\end{enumerate}
Model-theoretical notions will refer to $T$; if we mean them in the sense of $T_0$, we will indicate this by the index $0$. Moreover, we will assume that $T_0$ has geometric elimination of imaginaries, i.e.\ every $T_0$-imaginary element is $T_0$-interalgebraic with a real tuple. Note that this always holds if $T_0$ is strongly minimal with infinite $\acl_0(\emptyset)$.

It is easy (but nontrivial) to see that $A\ind_BC$ implies $A\ind^0_BC$ whenever $B$ is algebraically closed (in the sense of $T$).

In \cite{BMPW09} a relative version of CM-triviality was introduced: A theory $T$ is {\em CM-trivial over $T_0$ with respect to $\sscl.$} if for every real algebraically closed sets $A\subseteq B$ and every real tuple $\bar c$, whenever
\[ \sscl{A\bar c} \ind^0_{A} B,\]
the canonical base $\cb(\bar c/A)$ is algebraic over $\cb(\bar c/B)$ (in the sense of $T^\text{eq}$).

Every theory is CM-trivial over itself with respect to $\acl$. If $T$ is CM-trivial over its reduct to equality with respect to $\acl$, then $T$ is CM-trivial in the classical sense; the converse holds if $T$ has geometric elimination of imaginaries.

It follows from property (1) of the closure operator that relative CM-triviality is preserved under adding and forgetting parameters, and that we may assume $A$ and $B$ to be models.

A CM-trivial group of finite Morley rank is nilpotent-by-finite \cite{Pi95}; this remains true in a arbitrary stable group in the presence of enough regular types, or if the group is soluble \cite{Wa98}. For relative CM-triviality, a relative version was shown in \cite{BMPW09}:
\begin{quote} If $T$ is CM-trivial over $T_0$ with respect to $\sscl.$, then every connected type-definable group $G$ in $T$ allows a type-definable homomorphism to a group $H$ type-definable in $T_0$ whose kernel is contained (up to finite index) in the centre $Z(G)$ of $G$.\end{quote} 
It follows that a simple group or a field type-definable in $T$ embeds into one type-definable in $T_0$.
In fact, this consequence does not even need property (2).

\subsection{Relative CM-triviality of the coloured fields and the fusions}
If we consider the coloured fields or the fusions in a relational language with elimination of quantifiers, except possibly for the distinguished group law (addition for the red field, multiplication for the green field, vector space addition for the fusion over an $\FF_p$-vector space), the closure operator will satisfy conditions (1) and (2) above. Indeed, 
the amalgam of $B$ and $C$ over their algebraically closed intersection $A=B\cap C$ will be free if and only if
$$\begin{aligned}\mbox{{\bf (Fusion)}}\qquad&B\ind_A^{T_1} C\mbox{ and }B\ind_A^{T_2} C\\
\mbox{({\bf Red field)}}\qquad&B\ind_A^{ACF_p} C\mbox{ and }R(\langle BC\rangle)=R(B)+R(C)\\
\mbox{({\bf Green field)}}\qquad&B\ind_A^{ACF_0} C\mbox{ and }\U(\langle BC\rangle)=\U(B)\cdot\U(C).\end{aligned}$$
The characterization of independence in the generic model as
$$\bar b\ind_A\bar c\quad\Leftrightarrow\quad\parbox{6.5cm}{$\sscl{A\bar b}$ and $\sscl{A\bar c}$ are freely amalgamated over $A$ and the amalgam is strong}$$%
over any algebraically closed $A$ now implies property (1). Moreover, if $A$ is strong and $\bar a\in\acl_0(A)$, then $\delta(\bar a/A)\le 0$, so $A\bar a$ is strong. Let
$$B=\acl(\bar a)=\acl(\sscl{\bar a}).$$
Then $\cl(\bar a)\subseteq A\bar a\cap B$, whence $\delta(A\bar a\cap B)\ge\delta(\sscl{\bar a})=\delta(B)$,
and by submodularity
$$\delta(AB)\le\delta(A\bar a)+\delta(B)-\delta(A\bar a\cap B)\le\delta(A\bar a)\le\delta(A)~;$$
since $A$ is strong, so is $AB$. This yields property (2).

We shall want to check relative CM-triviality for the coloured fields (the proof for the fusion is analogous).

As for the absolute version there is an equivalent definition of relative CM-triviality, {\em non 2-ampleness}: For all real tuples $\bar a$, $\bar b$ and $\bar c$:
$$\acl(\bar a,\bar b) \ind^0_{\acl(\bar a)} \sscl{\acl(\bar a),\bar c}\quad\mbox{and}\quad
\bar a \ind_{\bar b}\bar c\quad\mbox{imply}\quad\bar a\ind_{\acl^{eq}(\bar a)\cap\acl^{eq}(\bar b)}\bar c.$$
So consider tuples $\bar a$, $\bar b$ and $\bar c$ such that:
\begin{enumerate}
\item $\bar a$ et $\bar b$ are algebraically closed,
\item $\acl(\bar a,\bar b) \ind^0_{\bar a} \sscl{\bar a,\bar c}$, and
\item $\bar a\ind_{\bar b}\bar c$.\end{enumerate}
Since relative CM-triviality is preserved under adding and forgetting parameters, we may add a $|T|^+$-saturated model $\M$ independent of $\bar a\bar b\bar c$ over $\acl^{eq}(\bar a)\cap\acl^{eq}(\bar b)$, and thus suppose $\acl^{eq}(\bar a)\cap\acl^{eq}(\bar b)=\acl^{eq}(\emptyset)$. We thus have to show $\bar a\ind\bar c$. We can also assume that $\bar a$ is a model.

Condition (3) means that $\acl(\bar a,\bar b)\ind^0_{\bar b}\acl(\bar b,\bar c)$ and that
$\langle\acl(\bar a,\bar b),\acl(\bar b,\bar c)\rangle$ is strong. We intersect this with the strong subset $\sscl{\bar a,\bar c}$, so the intersection is again strong.

Condition (2) implies $\sscl{\bar a,\bar b}\cap\acl(\bar a,\bar c)\subseteq\bar a$, whence
$$\sscl{\bar a,\bar c}\cap\bar b\subseteq \bar a\cap \bar b \subseteq \acl(\emptyset).$$
Put $D=\sscl{\bar a,\bar c}\cap\acl(\bar b,\bar c)\supset\sscl{\bar c}$. Then
$\cb_0(D/\sscl{\bar a,\bar b})\subseteq\bar a\cap\bar b=\acl(\emptyset)$.
Hence $\bar a\ind^0 D$.

In order to show that $\langle D\bar a\rangle$ is strong, consider $\bar\gamma\in\sscl{D\bar a}\setminus\langle D\bar a\rangle$ minimal with $\delta(\bar \gamma/D\bar a)<0$. As $\bar c\ind_{\bar b}\bar a$ is a free amalgam, there are $\bar\gamma_1\in\acl(\bar a,\bar b)$ and $\bar\gamma_2\in\acl(\bar c,\bar b)$ with $\bar\gamma=\bar\gamma_1\cdot\bar\gamma_2$.
Conditions (2) and (3) imply
$$D\bar\gamma\ind^0_{\bar a}\acl(\bar a\bar b)\quad\mbox{and}\quad D\bar\gamma_2\ind^0_{\bar b}\acl(\bar a\bar b).$$
Put $p_i(X,\bar x,\bar a)=\tp_i(D,\bar\gamma/\bar a)$ and let $E$ be the relation on $\tp(\bar a)$ given by
$$\bar a'E\bar a''\quad\Leftrightarrow\quad\exists\,\bar\gamma'\in\Gamma^{|\bar\gamma|}\ \bigwedge_{i<n}\bar\gamma'\cdot p_i(X,\bar x,\bar a')|_{\bar a',\bar\gamma'}\parallel p_i(X,\bar x,\bar a'')$$
(where $\parallel$ means parallelism of types: They have a---unique---common non-forking extension). This is a type-definable equivalence relation, and the class of $\bar a$ is definable over $\acl(\bar a)\cap\acl(\bar b)=\acl(\emptyset)$.

This enables us to find $\bar\gamma'_1\in\acl(\bar a)$ and $\bar\gamma'_2\in D$ with $\bar\gamma=\bar\gamma'_1\cdot\bar\gamma'_2$, so $\bar\gamma\in\langle D\bar a\rangle$, a contradiction. Thus $\langle D\bar a\rangle$ is strong.

If we choose $\bar\gamma\in\langle D\bar a\rangle$ to be coloured, we can do a similar argument and find first coloured $\bar\gamma_1\in\acl(\bar a,\bar b)$ and $\bar\gamma_2\in\acl(\bar c,\bar b)$, and then coloured $\bar\gamma_1'\in\acl(\bar a)$ and $\bar\gamma_2'\in D$ with $$\bar\gamma=\bar\gamma_1\cdot\bar\gamma_2=\bar\gamma_1'\cdot\bar\gamma_2'.$$
Thus $\bar a$ and $D$ are freely amalgamated, whence $\bar a\ind D$ and $\bar a\ind \bar c$.

It follows that in a coloured field every simple definable group is linear; in the fusion of strongly minimal $T_1$ and $T_2$ (over an $\aleph_0$-categorical reduct) every simple group is $T_1$-definable or $T_2$-definable. More generally, in a coloured field every definable group embeds modulo a central subgroup into an algebraic group; in the fusion of strongly minimal $T_1$ and $T_2$ every definable group embeds modulo a central subgroup into a product of a $T_1$-definable and a $T_2$-definable group.

\subsection{Subgroups of groups definable in a reduct}
It remains to be seen that a simple $T$-definable subgroup $H$ of a $T_0$-definable group $G$ is actually $T_0$-definable, both in the fusion and in a coloured field. For a coloured field, let $a$ and $b$ be two independent generic elements of $H$, and $A=\cl(a,\acl(\emptyset))$ and $B=\cl(b,\acl(\emptyset))$. Since $A\ind B$, they are freely amalgamated over $A\cap B=\acl(\emptyset)$ and $\langle AB\rangle$ is strong. Let $c=ab$, and $C=\cl(c,\acl(\emptyset)$. Since $c\in\acl_0(a,b)$ we have $\delta(c/\langle AB\rangle)\le 0\,$; as $\langle AB\rangle$ is strong, $\langle ABc\rangle$ is strong and contains $C$. But $\tr(c/AB)=0$, whence $\lin(c/AB)=0$ and all coloured points $C_0$ of $C$ are already in $\langle AB\rangle$, where they are the sum (or product) of coloured points $A_0$ of $A$ with coloured points $B_0$ of $B$. Since $A,B,C$ are pairwise independent, so are $A_0,B_0,C_0$; as $A\equiv B\equiv C$ we must have $A_0\equiv B_0\equiv C_0$. Ziegler's Lemma now implies that $A_0$, $B_0$ and $C_0$ realize the generic type of a $T$-definable connected subgroup $V$ of some Cartesian power of $R$ (or $\U$). But then the correspondence $a\mapsto A_0$ induces a $T$-definable homomorphism $\phi:H\to V$ (modulo a finite co-kernel), which must be trivial by simplicity of $H$. Hence $A$ has no coloured points outside of $\acl(\emptyset)$; quantifier-elimination now implies that $H$ is $T_0$-definable.

More generally, the argument shows that a connected subgroup $H$ of an algebraic group $G$ definable in a coloured field is an extension of a subgroup $V$ of some Cartesion power of the colour subgroup by an algebraic subgroup $N=\ker\phi$ of $H$; in the collapsed case or in the green field $V$ must itself be a Cartesian power of the colour subgroup, by strong minimality of $R$ or degeneracy of pre-algebraic extensions in the green field, respectively.

In the fusion of strongly minimal $T_1$ and $T_2$ over a common $\FF_p$-vector space $V$, let $H$ be a $T$-definable connected subgroup of a $T_1$-definable group $G$, and choose $a,b,c,A,B,C$ as above. Put $C_0=C\cap\langle AB\rangle$. Since $\langle AB\rangle$ is a free amalgam, there is $A_0\subseteq A$ and $B_0\subseteq B$ with $C_0=A_0+B_0$ (coordinatewise); again $A\equiv B\equiv C$ implies $A_0\equiv B_0\equiv C_0$. By Ziegler's Lemma $A_0$, $B_0$ and $C_0$ realize the generic type of a $T$-definable connected subspace $V_0$ of some Cartesian power of $V$, and $a\mapsto A_0$ induces a homomorphism $\phi:H\to V_0$ with kernel $N=\ker\phi$.  Now $RM_1(c/AB)=0$, so $\langle ABc\rangle$ is strong and contains $C$, whence $RM_2(C/AB)=RM_0(C/AB)$. But $\tp(c/C_0)$ is the generic type of the coset $cN$; since
$$RM_0(C/C_0)\ge RM_2(C/C_0)\ge RM_2(C/AB)=RM_0(C/AB)=RM_0(C/C_0)$$
(by modularity of $RM_0$), equality holds everywhere and $\tp(C/C_0)$ is $T_2$-free, i.e.\ generic in $\tp_1(C/C_0)$. This implies that $N$ is $T_1$-definable, and $G$ is an extension of $V_0$ by $N$. If $H$ is simple, $N$ must be trivial, so $H$ itself is $T_1$-definable. In the collapsed case, $V_0$ is isomorphic to a Cartesion power of $V$ by strong minimality.

The argument can be adapted for a $T$-definable connected subgroup $H$ of a product $G=G_1\times G_2$, where $G_i$ is $T_i$-definable for $i=1,2$. Again we take $a=(a_1,a_2)$ and $b=(b_1,b_2)$ generic points of $H$, and we put $ab=c=(c_1,c_2)=(a_1b_1,a_2b_2)$. We define $A,B,C,A_0,B_0,C_0,V_0,\phi,N$ as in the previous paragraph, and put $C_i=\cl(c_i,\acl(C_0))$ for $i=1,2$. Note that $\langle AB\acl(C_0)\rangle$ is strong by property (2) of the closure operator.
$$RM_1(c_1/AB\acl(C_0))=RM_2(c_2/AB\acl(C_0))=0,$$
so $\langle ABc_1\acl(C_0)\rangle$ and $\langle ABc_2\acl(C_0)\rangle$ are both strong and contain $C_1$ and $C_2$, respectively.
Hence
$$RM_2(C_1/AB\acl(C_0))=RM_0(C_1/AB\acl(C_0))$$
and
$$RM_1(C_2/AB\acl(C_0))=RM_0(C_2/AB\acl(C_0)),$$
whence
$$\begin{aligned}RM_0(C_1/\acl(C_0))&\ge RM_2(C_1/\acl(C_0))\ge RM_2(C_1/C_2)\\
&\ge RM_2(C_1/ABC_2)=RM_2(C_1/AB\acl(C_0))\\
&=RM_0(C_1/AB\acl(C_0))=RM_0(C_1/\acl(C_0)).\end{aligned}$$
It follows that
$$C_1\ind^2_{\acl(C_0)}C_2\quad\mbox{and, similarly,}\quad C_2\ind^1_{\acl(C_0)}C_1.$$
So $C_1$ and $C_2$ are freely amalgamated over $\acl(C_0)$.
Moreover, $\langle AB\acl(C_0)c_1c_2\rangle$ and $\langle C\acl(C_0)\rangle$ are both strong; intersecting them we obtain that $\langle C_1C_2\rangle$ is strong. Thus $C_1$ and $C_2$ are independent over $\acl(C_0)$. This means that $N$ is the product of its projections to $G_1$ and $G_2$, and thus equal to the product of some $T_1$-definable $N_1\le G_1$ and some $T_2$-definable $N_2\le G_2$.

{\small\sc Universit\'e de Lyon\\
Universit\'e Lyon 1\\
CNRS\\
Institut Camille Jordan UMR 5208\\
Villeurbanne\\
F-69622\\
France}

\begin{thebibliography}{99}
\bibitem{BH} J. Baldwin, K. Holland,
\emph{Constructing $\omega$-stable structures: 
fields of rank 2}, J. Symb.\ Logic, {\bf 65} (2000), 371--391.
\bibitem{BH2} J. Baldwin, K. Holland,
\emph{Constructing $\omega$-stable structures: 
Computing Rank}, Fund.\ Math., {\bf 170} (2001), 1--20.
\bibitem{BH3} J. Baldwin, K. Holland,
\emph{Constructing $\omega$-stable structures: Model
Completeness}, Ann.\ Pure Appl.\ Logic, {\bf 125} (2004), 159--172.
\bibitem{BS} J. Baldwin, S. Shelah, \emph{Randomness and semigenericity},
Trans.\ Amer.\ Math.\ Soc., {\bf 349} (1997), 1359--1376.
\bibitem{Bau} A. Baudisch, \emph{A new uncountably categorical group}, Trans.\
Amer.\ Math.\ Soc., {\bf 348} (1996), 3889--3940. 
\bibitem{BHMW} A. Baudisch, M. Hils, A. Mart\'\i n Pizarro, F.
O. Wagner, \emph{Die {B}\"{o}se {F}arbe}, J. Inst.\ Math.\ Jussieu, {\bf 8} (2009), 415--443.
\bibitem{BMPZfields} A. Baudisch,
A. Mart{\'\i}n Pizarro, M. Ziegler, \emph{On fields and colours}, Algebra i
Logika,  {\bf 45} (2006), 92--105.
\bibitem{BMPZredfus} A. Baudisch, A. Mart{\'\i}n Pizarro,
M. Ziegler, \emph{Fusion over a vector space}, J. Math.\ Logic, {\bf 6} (2006), 141--162.
\bibitem{BMPZfus} A. Baudisch, A. Mart{\'\i}n Pizarro, M. Ziegler, \emph{Hrushovski's Fusion}. In:
F. Haug et al (eds), Festschrift f{\"u}r Ulrich Felgner zum 65.\ Geburtstag, pp.\ 15--31, Studies in Logic 4, College Publications, 2007.
\bibitem{BMPZred} A. Baudisch, A. Mart{\'\i}n Pizarro, M. Ziegler,
\emph{Red fields}, J. Symb.\ Logic, {\bf 72} (2007), 207--225.
\bibitem{BMPW09} T. Blossier, A. Mart{\'\i}n Pizarro, F. Wagner, \emph{G\'eom\'etries relatives}, preprint.
\bibitem{BMPW10} T. Blossier, A. Mart{\'\i}n Pizarro, F. Wagner, \emph{Groupes d\'efinissables dans les fusions de Hrushovski et dans les corps color\'es}, in preparation.
\bibitem{Ev1} D. Evans, \emph{Ample dividing}, J. Symb.\ Logic, {\bf
68} (2003), 1385--1402.
\bibitem{Ev2} D. Evans, \emph{Trivial stable structures with non-trivial
reducts}, J. Lond.\ Math.\ Soc.\ (2), {\bf 72} (2005), 351--363. 
\bibitem{Ha1} A. Hasson, \emph{In search of new strongly minimal sets},
Ph.D. thesis, Jerusalem (2004). 
\bibitem{Ha3} A. Hasson, \emph{Interpreting structures of finite Morley rank in strongly minimal sets}, Ann.\ Pure Appl.\ Logic {\bf 145} (2007), 96--114.
\bibitem{Ha} A. Hasson, \emph{Some questions concerning Hrushovski's amalgamation constructions}, J. Inst.\ Math.\ Jussieu {\bf 7} (2008), 793--823. 
\bibitem{Ha2} A. Hasson, \emph{Collapsing structures and a theory of
envelopes}, preprint. 
\bibitem{HH} A. Hasson, M. Hils, \emph{Fusion over sublanguages},
J. Symb.\ Logic {\bf 71} (2006), 361--398.
\bibitem{Hi} M. Hils, \emph{Fusion libre et autres constructions
g\'en\'eriques}, Ph.D. thesis, Universit\'e Paris-7, Paris 2006. 
\bibitem{Hi1} M. Hils, \emph{La fusion libre: le cas simple}, J. Inst.\ Math. Jussieu {\bf 7} (2008), 825--868.
\bibitem{Hi2} M. Hils, \emph{Generic automorphisms and green fields}, in preparation.
\bibitem{Hr92} E. Hrushovski, \emph{Strongly minimal expansions of
algebraically closed fields}, Isr.\ J. Math., {\bf 79} (1992), 129--151.
\bibitem{Hr93} E. Hrushovski, \emph{A new strongly minimal
set}, Ann.\ Pure Appl.\ Logic, {\bf 62} (1993), 147--166.
\bibitem{Ik} K. Ikeda, \emph{Minimal but not strongly minimal structures with arbitrary finite dimensions},
J. Symb.\ Logic {\bf 66} (2001), 117--126. 
\bibitem{I} K. Ikeda, \emph{The stability spectrum of ab initio generic structures}, preprint 2010.
\bibitem{L} C. Laskowski, \emph{A simpler axiomatization of the Shelah-Spencer
almost sure theory}, Isr.\ J. Math., {\bf 161} (2007), 157--186.
\bibitem{Pi95} A. Pillay, \emph{The geometry of forking and groups of finite Morley rank}, J. Symb.\ Logic {\bf 60} (1995), 1251--1259.
\bibitem{Po99} B. Poizat, \emph{Le carr\'e
de l'\'egalit\'e }, J. Symb.\ Logic, {\bf 64} (1999), 1338--1355.
\bibitem{Po01} B. Poizat, \emph{L'\'egalit\'e au
cube}, J. Symb.\ Logic, {\bf 66} (2001), 1647--1676.
\bibitem{ShSp} S. Shelah, J. Spencer, \emph{Zero-one laws for sparse random graphs},
J. Amer.\ Math.\ Soc., {\bf 1} (1988), 97--115.
\bibitem{Wa98} F. O. Wagner, \emph{{CM}-triviality and stable groups}, J. Symb.\ Logic {\bf 63} (1998), 1473--1495.
\bibitem{Wa03} F. O. Wagner, \emph{Bad fields in positive characteristic}, Bull.\ London Math.\ Soc., {\bf 35} r(2003), 499--502.
\bibitem{wa06} F. O. Wagner,\emph{Hrushovski's amalgamation construction}. In: S. B. Cooper et al (eds.), Logic Colloquium 2006, pp.\ 361-373, Ass.\ Symb.\ Logic LN Logic 32, Cambridge University Press, 2009.
\bibitem{Zi} M. Ziegler, \emph{Fusion of structures of finite Morley
rank}. In: Z. Chatzidakis et al (eds.),  Model theory with Applications to Algebra and Analysis Vol.\ 1, pp.\ 225--248, London Math.\ Soc.\ LN 350, Cambridge University Press, 2008.
\bibitem{ZiNote} M. Ziegler, \emph{A Note on generic Types} (2006). http://arxiv.org/math.LO/0608433.
\end{thebibliography}
\end{document}